\documentclass[11pt]{article}

\input{amssym.def}
\input{amssym}
\usepackage{eufrak}
\usepackage{color}

\newtheorem{theorem}{Theorem}

\newtheorem{proposition}[theorem]{Proposition}
\newtheorem{remark}[theorem]{Remark}

\def\RRR{{\cal{R}}}

\def\Y{\mathbf{Y}}
\def\YR{\mathbf{Y}(R)}

\def\YYR{\mathbf{Y}_{RRT}(R)}

\def\qq{q^{-1}}

\def\lhq{{\cal L}(\h, q)}

\def\h{\hbar}
\def\lrh{{\cal L}(R,h)}
\def\lro{{\cal L}(R,1)}
\def\lr{{\cal L}(R)}

\def\TT{{\mathbf T}}
\def\LL{{\mathbf L}}
\def\T{{\cal T}(R)}

\def\de{\delta}

\def\gg{\mbox{$\frak g$}}
\def\ot{\otimes}

\def\K{{\Bbb K}}

\def\C{{\Bbb C}}
\def\R{{\Bbb R}}

\def\vv{V^{\otimes 2}}

\def\lhq{\ifmmode {\cal L}(q,\hbar)\else ${\cal L}(q,\hbar)$\fi}

\def\lqh{\ifmmode {\cal L}(q,\hbar)\else ${\cal L}(q,\hbar)$\fi}

\def\al{{\alpha}}

\def\be{\begin{equation}}
\def\ee{\end{equation}}

\begin{document}

\makeatletter
\renewcommand{\theequation}{{\thesection}.{\arabic{equation}}}
\@addtoreset{equation}{section} \makeatother

\title{Generalized Yangians and their Poisson counterparts}

\author{\rule{0pt}{7mm} Dimitri
Gurevich\thanks{gurevich@ihes.fr}\\
{\small\it LAMAV, Universit\'e de Valenciennes,
59313 Valenciennes, France}\\
\rule{0pt}{7mm} Pavel Saponov\thanks{Pavel.Saponov@ihep.ru}\\
{\small\it
National Research University Higher School of Economics,}\\
{\small\it 20 Myasnitskaya Ulitsa, Moscow 101000, Russian Federation}\\
{\small \it and}\\
{\small \it
Institute for High Energy Physics, NRC "Kurchatov Institute"}\\
{\small \it Protvino 142281, Russian Federation}}

\maketitle
\vskip 5mm

\rightline{\it To the memory of our friend Peter Kulish}
\vskip 5mm

\begin{abstract}
By a generalized Yangian we mean a Yangian-like algebra of one of two classes. One of these classes consists of the so-called braided Yangians, introduced in our previous paper.
The braided Yangians are in a sense similar to the reflection equation algebra. The generalized Yangians of second class, called the Yangians of RTT type, are defined by the same
formulae as the usual Yangians are but with other quantum $R$-matrices. If such an $R$-matrix is the simplest trigonometrical $R$-matrix, the corresponding Yangian
of RTT type is the so-called q-Yangian. We claim that each generalized Yangian is a deformation of the commutative algebra ${\rm Sym }(gl(m)[t^{-1}])$ provided that the corresponding
$R$-matrix is a deformation of the flip. Also, we exhibit the corresponding Poisson brackets.
\end{abstract}

{\bf AMS Mathematics Subject Classification, 2010:} 81R50


\section{Introduction}

In \cite{GS} we introduced the notion of the braided Yangians associated with a wide class of rational and trigonometrical $R$-matrices. This notion is a new generalization of the
 Yangian $\Y(gl(m))$ introduced by Drinfeld \cite{D1}. According to one of the equivalent definitions, Drinfeld's Yangian $\Y(gl(m))$
is the algebra generated by the coefficients of the matrix-valued function
 \be L(u)=\sum_{k\geq 0} L[k]u^{-k}, \label{L}
\ee
subject to the system
 \be
R(u,v)L_1(u)L_2(v)-L_1(v)L_2(u)R(u,v)=0,
\label{sys}
\ee
where $R(u,v)=P-\frac{a}{u-v}I$ is the famous Yang $R$-matrix\footnote{Note that this $R$-matrix is often used in a different form ${\cal R}(u,v) = PR(u,v)$.}.
Hereafter, $I$ is the identity matrix, $P$ is the usual flip or its matrix, and
the Laurent coefficients $L[k]=\|l_i^j[k]\|$ are  matrices. Also, a complementary condition $L[0]=I$ is imposed.
All matrices are of size $m\times m$. The matrix $L(u)$ (and all similar matrices, considered below)
will be called {\em generating}.

The Yang $R$-matrix entering this definition is the simplest example of current (i.e. depending on parameters) $R$-matrices, which are, by definition, solutions to
 the so-called quantum Yang-Baxter equation

\be R_{12}(u, v)R_{23}(u, w)R_{12}(v, w)=R_{23}(v, w)R_{12}(u, w)R_{23}(u, v). \label{YB} \ee

If in the relation (\ref{sys}) the Yang $R$-matrix is replaced by another current $R$-matrix $R(u,v)$, we get a class of Yangian-like algebras. The best known example is the so-called
{\em q-Yangian}, which corresponds to the simplest trigonometrical $R$-matrix related to the quantum group (QG) $U_q(\widehat{sl(m)})$. Below, all objects, related to $U_q(\widehat{sl(m)})$ or $U_q(sl(m))$ are called {\em standard}.

However, there exists a large family of rational and trigonometrical current $R$-matrices, which can
be constructed from involutive or Hecke symmetries $R$ by means of the Yang-Baxterization procedure (see Section 2 for the definitions).
We call the corresponding algebras (\ref{sys})  {\em Yangians of RTT type} and denote them $\YYR$. Note that by contrast with the usual Yangian, the condition $L[0]=I$ is not imposed.

Another class of algebras --- the braided Yangians --- was introduced in \cite{GS}. Each of them is defined by the following system
\be
R(u,v)L_1(u)R L_1(v)-L_1(v)R L_1(u)R(u,v)=0,
\label{sys1}
\ee
where $R(u,v)$ is the same current $R$-matrix resulting from the Yang-Baxterization of the symmetry $R$.
Here, the matrix $L(u)$ is also a series given by the formula (\ref{L}). Besides, the condition $L[0]=I$ is imposed.
The braided Yangian arising from a symmetry $R$ (involutive or Hecke) will be denoted $\YR$.

The above Yangian-like algebras of both classes are called {\em generalized Yangians}. Below we often omit the term {\em generalized}.

Note that the properties of these two classes of the generalized Yangians are distinct.
In particular, they have different bi-algebra structures. While the bi-algebra structure of a Yangian of RTT type is usual, that of a braided Yangian is more complicated:
the product and coproduct in such an algebra are coordinated via the corresponding braiding $R$, similarly to the case of super-algebras, where the role of $R$ is played by a super-flip.
Also, the evaluation maps and the corresponding target algebras differ drastically.

In a sense, the evaluation maps for the braided Yangians are more similar to the classical case. For a braided Yangian $\YR$ the corresponding target algebra is the reflection equation
(RE) algebra (may be, in its {\it modified} form) related to the initial symmetry $R$. Treating the modified RE algebra as a "braided analog" of the algebra $U(gl(m))$, we succeeded in
constructing its representation category, similar to that of $U(gl(m))$ (see \cite{GPS} for details). Using such categories and evaluation maps, we obtain large representation categories
of the corresponding braided Yangians.

For the Yangians of RTT type the target algebras are less interesting and their representation theories are not known. We refer the reader to \cite{Mo}, where the q-Yangian, i.e. the
Yangian of RTT type related to the standard trigonometrical $R$-matrix is considered as an example.

The main purpose of this paper is to study the {\em deformation property} of the generalized
Yangians of both classes\footnote{Note that all algebras which we are dealing with are defined via
generators and relations and the parameters coming in the relations are not formal and can be specialized.}. We say that an associative algebra $A_h$ depending on a deformation
parameter $h$ has the deformation property, if for $h=0$ it turns into a commutative algebra $A=A_0$ and in the latter algebra a new product $\star_h$ induced from the algebra $A_h$
and smoothly depending on $h$ can be defined. This means that for a generic $h$ there exists a linear space map $\al_h:A\to A_h$ smoothly depending on $h$ and such that
$\alpha_0 = {\rm Id}$. Then the induced product in the algebra $A$ is
$$
f\star_h g=\al_h^{-1}(\al_h(f)\circ \al_h(g)),
$$
where $\circ$ stands for the product in the algebra $A_h$. Such a map $\al_h$ is usually constructed by means of a special basis (sometimes called a PBW basis) in the algebra $A_h$.
Below we claim that each generalized Yangian is a deformation of the commutative algebra ${\rm Sym }(gl(m)[t^{-1}])$ provided that the symmetry $R$ is a deformation of the flip.

If an algebra $A_h$ has the deformation property, the product $\star_h$ can be expended in a series in $h$:
$$
f\star_h g= f\cdot g+h \, c_1(f, g)+h^2\,c_2(f, g)+\dots,
$$
where $\cdot$ stands for the commutative product in the algebra $A$. Then there exists the corresponding Poisson bracket defined by the skew-symmetrized term $c_1$:
$$
A^{\ot 2}\ni f\ot g \mapsto \{f,g\}=\frac{1}{2}(c_1(f,g)-c_1(g,f))\in A.
$$

The second purpose of the paper is to exhibit the  Poisson structures corresponding to the Yangians of both classes.
Namely, we compute the quadratic Poisson brackets arising from these Yangians
 and their linear counterparts, obtained by the linearization procedure.

Our paper is organized as follows. In the next section we recall the definitions of the RTT and RE algebras and
Yangians of both classes and compare some their properties. Section 3 is devoted to the study of the deformation property of the generalized Yangians.
Besides, the corresponding quadratic Poisson brackets are calculated.
Their linear counterparts are exhibited
in the last section. Also, two examples are given. The ground field $\K$ is assumed to be $\C$ or $\R$.

\medskip

\noindent
{\bf Acknowledgements.} D.G. is grateful to SISSA (Italy), where  during his scientific visit this paper was conceived. He is also grateful to V.Rubtsov for elucidating discussions.

P.S. is grateful to IHES where this work was completed during his scientific visit. The work of P.S. has been funded by the Russian Academic Excellence Project '5-100' and was 
also partially supported by the RFBR grant 16-01-00562.

\section{Quantum matrix algebras and generalized Yangians}

First, we consider quantum matrix (QM) algebras similar to the Yangians of both classes but associated with {\em constant braidings}. Recall that
an operator $R :\vv\to \vv$, where $V$ is a finite dimensional vector space ($\dim V = m$) is called a {\em braiding} if it meets the quantum Yang-Baxter equation (\ref{YB}) but without parameters.
A braiding $R$ is called {\em Hecke (resp., involutive) symmetry}, if it meets the relation
\be
(qI-R)(\qq I+R)=0,
\label{H-cond}
\ee
where the numeric nonzero parameter $q$ obeys the condition $q^2\not=1$ (resp., $q^2=1$). Below, if $q\not=1$, $q$ is assumed to be {\em generic}, that is $q^n\not=1$ for
any integer $n$.

The algebra, generated by entries of a matrix $L=\|l_i^j\|_{1\leq i, j\leq m}$ subject to the system
\be
R L_1 R L_1-L_1 R L_1 R=h(R L_1-L_1 R),
\label{mRE}
\ee
is called the {\em modified reflection equation} algebra if $h\not=0$. If $h=0$, we omit the term {\em modified}. The algebra defined by (\ref{mRE}) is denoted $\lrh$ if $h\not=0$ and $\lr$ if $h=0$.
As usual, the low indexes (in $L_i$, $R_{i\, j}$ and so on) indicate the positions in $V^{\ot n}$, where a given operator or matrix is located. The operator $R_{i\,i+1}$ will be also denoted $R_i$.

Besides, we consider the RTT algebra associated with a braiding $R$ and defined by the following system
\be
R T_1 T_2-T_1 T_2 R=0, \qquad T=\|t_i^j\|_{1\leq i, j\leq m}.
\label{RTT}
\ee
This algebra will be denoted $\T$.

Note that the algebras $\lr$ and $\T$ are particular cases of the QM algebras defined in general through two
braidings in a sense coordinated (see \cite{IOP}).

A braiding $R$ is called {\em skew-invertible} if there exits an operator $\Psi:\vv\to \vv$ such that
\be
{\rm Tr}_2 R_{12} \Psi_{23}=P_{13}\quad\Leftrightarrow\quad R_{ij}^{kl} \Psi_{lp}^{jq}=\de_i^q \de_p^k.
\label{def:Psi}
\ee
All braidings $R$ we are dealing with are assumed to be skew-invertible Hecke or involutive symmetries.

Emphasize that if $R$ is a skew-invertible braiding, it is possible to define the so-called $R$-trace ${\rm Tr}_R A$, where $A$ is an arbitrary $m\times m$ matrix.
This trace possesses a lot of remarkable properties. In particular, the $R$-trace enters a construction
of {\em quantum symmetric polynomials} in the algebras $\T$ and $\lr$ as well as in the generalized Yangians of both classes. Their explicit form
will be presented below for the generalized Yangians. Note, that these quantum symmetric polynomials generate
the so-called characteristic subalgebras $Ch(\T)$ and $Ch(\lr)$ of the algebras $\T$ and $\lr$ respectively. However, properties of the algebras $Ch(\T)$ and $Ch(\lr)$
are different: if the latter algebras are central,
the former algebras are not central, but only commutative. The reader is referred to \cite{IOP, IP} for details.

Let us turn, now, to the generalized Yangians. To this end we first specify the current $R$-matrices $R(u,v)$ entering formulae
(\ref{sys}) and (\ref{sys1}). In \cite{GS} the following claim was proved.

\begin{proposition} Consider the sum\rm
\be
R(u,v)=R+g(u,v)I,
\label{Rmat}
\ee
\it
where $R$ is a braiding, $g(u,v)=f(u-v)$, and the $f(z)$ is a non-constant
meromorphic function. If $R$ is an involutive symmetry, then $R(u,v)$ is a current $R$-matrix iff\rm
\be
g(u,v)=\frac{a}{u-v}.
\label{g}
\ee
\it
If $R=R_q$ is a Hecke symmetry, then $R(u,v)$ is a current $R$-matrix iff\rm
\be
 g(u,v)=\frac{q-\qq}{b^{u-v}-1}.
\label{gg}
\ee
\it
Here $a$ and $b\not=1$ are arbitrary nonzero complex numbers.
\end{proposition}

In particular, by putting $b=q^{-2/a}$ we get
\be
R(u,v)=R_q-\frac{q^{\frac{u-v}{a}}}{(\frac{u-v}{a})_q}I.
\label{trig}
\ee
As $q\to 1$ this $R$-matrix tends to
\be
R_1-\frac{a}{u-v}\, I,
\label{rrat}
\ee
provided the Hecke symmetry $R_q$ tends to an involutive one $R_1$.

Changing the variables $b^{-u}\to u$, $b^{-v}\to v$ in (\ref{gg}), we get the following form of the trigonometric current $R$-matrix
\be
R(u,v)=R-\frac{u (q-\qq)}{u-v}I.
\label{rat}
\ee
Evidently, it depends only on the ratio $x=v/u$.

Below, we deal with the Yangians of RTT type $\YYR$ and braided Yangians $\YR$, respectively defined by formulae (\ref{sys}) and (\ref{sys1}), where the factors $R(u,v)$ stand for the
current $R$-matrices (\ref{rrat}) or (\ref{rat}) and the middle terms $R$ in (\ref{sys1}) are the initial symmetries.

As we noticed above, each of the Yangians of both classes has a bi-algebra structure, but in $\YR$ this structure is braided. We refer the reader to \cite{GS} for details.
Here, we exhibit the evaluation morphisms for the braided Yangians. For a given braided Yangian $\YR$ the evaluation morphism is defined by the following formula
\be L(u) \mapsto I+\frac{M}{u}, \label{ev} \ee
where $M$ is the generating matrix of the {\em target algebra}. A concrete form of this algebra depends on
the initial symmetry $R$.

\begin{proposition} {\bf \cite{GS}}
\begin{enumerate}
\item[\rm 1.] If $R$ is an involutive symmetry, then the map (\ref{ev})
defines a surjective morphism $\YR\to \lro$. Besides, the map $M\mapsto L[1]$ defines an injective morphism $\lro\to \YR$.

\item[\rm 2.] If $R$ is a Hecke symmetry, then the map (\ref{ev})
defines a morphism $\YR\to \lr$.
\end{enumerate}
Thus, the target algebra is $\lro$ in the former case and $\lr$ in the latter case.
\end{proposition}

This proposition enables us to construct a category of finite dimensional representations of a braided Yangian by using the results of the paper \cite{GPS}, where such a category was
 constructed for the algebra $\lro$. Thus, if $R$ is an involutive symmetry, the evaluation map converts any
$\lro$-module into a $\YR$-one. If $R$ is a Hecke symmetry, then   first  we convert any $\lro$-module into $\lr$-one. It is possible to do since the algebras $\lr$ and
$\lro$ are isomorphic to each other. Their isomorphism is given by the following map
$$
\lr \to \lro,\quad L\mapsto L-\frac{1}{q-\qq} I,\quad q^2\not=1.
$$

By contrast, the evaluation morphism for the Yangians of RTT type has the following form
\be L(u)\mapsto T_0+\frac{T_1}{u}. \label{ev1}\ee
If $R$ is a Hecke symmetry,
the map (\ref{ev1}) leads to a target algebra, which is similar to the quantum algebra $U_q(gl(m))$. The explicit form of this quantum algebra
 can be found, for instance, in \cite{Mo,FJMR}. If $R$ is an involutive symmetry,
the target algebra is new and its properties are not known. (Note that the condition $T_0=I$ leads to an inconsistent target algebra. This is the reason, why the condition $T[0]=I$ is not imposed while defining the Yangians of RTT type.)

Now, we describe the aforementioned {\em quantum symmetric polynomials} in the braided Yangians. Also, we write down the Newton-Cayley-Hamilton (NCH) identities in these algebras
in the spirit of \cite{IOP}. To this end we first introduce quantum analogs of matrix powers and skew-powers. Here, we deal with
 trigonometrical $R$-matrices (\ref{trig}) where we put $a=1$. The case of rational $R$-matrices (\ref{rrat}) can be obtained after the passage to the form (\ref{trig}) of the $R$-matrix.

Let us define the {\em quantum skew-powers} of the matrix $L(u)$ as follows
 \be
L^{\wedge k}(x) = {\rm Tr}_{R(2\dots k)}\Big({\cal P}^{(k)}_{12\dots k}L_{\overline 1}(x)L_{\overline 2}(x-1)\dots L_{\overline k}(x-k+1)\Big),
\quad k\ge 2,
\label{wedge-p}
\ee
where
\begin{eqnarray}
{\cal P}^{(k+1)}_{12\dots k+1}&=&\frac{(-1)^k}{(k+1)_q}\,R_1(1)R_2(2)\dots R_k(k)\,{\cal P}^{(k)}_{12\dots k}
\end{eqnarray}
is the skew-symmetrizer in the space $V^{\otimes (k+1)}$ associated with the Hecke symmetry $R$. Also, we put by definition $L^{\wedge 1}(x) = L(x)$.

Note, that if the Hecke symmetry $R$ is a deformation of the usual flip $P$, then $L^{\wedge k}(x)\equiv 0$ for
$k>m$.

Also, define the {\em quantum matrix powers} of the generating matrix as follows
\be
L^{k}(x) = {\rm Tr}_{R(2\dots k)}\Big(R_1R_2\dots R_{k-1}L_{\overline 1}(x)L_{\overline 2}(x-1)\dots L_{\overline k}(x-k+1)\Big),
\quad k\ge 1.
\label{power-p}
\ee

Now, the {\em quantum elementary symmetric polynomials} and {\em quantum power sums} are respectively defined by
\be
e_k(x) = {\rm Tr}_R(L^{\wedge k}(x)),\qquad s_k(x) = {\rm Tr }_R(L^k(x)).
\label{pow-sum}
\ee
\begin{proposition}\label{prop:22}
The quantum skew-symmetric powers and matrix powers of $L(u)$ are related by the Cayley-Hamilton-Newton identities \rm
\be
(-1)^{k+1}k_qL^{\wedge k}(x) = \sum_{p=1}^k(-q)^{k-p}L^p(x)e_{k-p}(x-p),\quad k\ge 1.
\label{chn-id}
\ee
\end{proposition}

If $R$ is a deformation of $P$, then the last non-trivial {\em Cayley-Hamilton-Newton identity} turns into the quantum {\em Cayley-Hamilton identity}.
In this case the highest nontrivial symmetric polynomial $e_m(x)$ is called the {\em quantum determinant}.

Applying the $R$-trace to the relation (\ref{chn-id}), we get relations between the quantum powers sums and quantum elementary symmetric polynomials.

Note that in the Yangians of RTT type all these objects (quantum matrix powers, skew-powers, power sums and elementary symmetric polynomials) are also well defined. Observe that
similarly to the RTT algebras the quantum elementary symmetric polynomials and power sums in the Yangians of both classes generate commutative subalgebras.
This fact will be proved in our subsequent publications.

 Nevertheless, by contrast with the RE algebras, these commutative subalgebras are not central.
The only non-trivial quantum symmetric polynomial, which is central in the braided Yangian $\YR$, is the quantum determinant. To be more precise, all coefficients of its expansion
in a series are central. As for the Yangian $\YYR$, its quantum determinant
is central iff it is so in the corresponding RTT algebra.
Thus, the centrality of the quantum determinant in the Yangian of RTT type
depends on the initial symmetry $R$.

In the last section we give an example of an RTT algebra and consequently
of the Yangian $\YYR$ with non-central quantum determinant.

\section{Deformation property of Yangians and corresponding Poisson structures}

Since all Yangians we are dealing with are introduced via generators and relations, we have to check their deformation property, provided that a symmetry $R$ (involutive or Hecke)
is a deformation of the flip. Below we present arguments in favour of the deformation property of the generalized Yangians. Also, we exhibit  the corresponding Poisson structures.

However, first we consider the deformation property of the QM algebras $\T$ and $\lr$. To this end we introduce the following notations
$$
L_{\overline 1}=L_1,\qquad L_{\overline k}=R_{k-1 } L_{\overline{k-1}} R^{-1}_{k-1},\quad \quad \forall\,k\ge 2.
$$
Then for $h=0$ the relations (\ref{mRE}) can be written in a form similar to the relations in an RTT algebra:
$$R L_{\overline 1} L_{\overline 2}=L_{\overline 1} L_{\overline 2} R.$$
These notations are also useful for establishing some linear space isomorphism between the algebras $\T$ and $\lr$.
Respectively, denote
$$
\TT={\rm span}_{\K}(t_i^j)\qquad \LL={\rm span}_{\K}(l_i^j)
$$
the vector spaces spanned by the generators $t_i^j$ of the RTT algebra $\T$ and by the generators $l_i^j$ of the RE algebra $\lr$, associated with a symmetry $R$.
For any positive integer $k$ we consider the linear map $\pi_k: \TT^{\ot k}\to \LL^{\ot k}$ of the corresponding tensor powers which is defined on the basis elements
by the following rule:
$$
\pi_k(T_1\ot T_2\ot ...\ot T_k) =L_{\overline 1}\ot L_{\overline 2}\ot \dots\ot L_{\overline k}, \quad k\ge 1.
$$
Below, we omit the tensor product sign $\otimes$ to simplify the formulae.
\begin{proposition} The following relations take place:\rm
\begin{eqnarray}
\pi_k (T_1\dots T_{i-1} (R_i T_i T_{i+1}- T_i T_{i+1} R_i )T_{i+1}\dots T_k) &&\nonumber\\
&&\hspace*{-30mm} =L_{\overline 1}\dots L_{\overline {i-1}} (R_i L_{\overline i} L_{\overline {i+1}}-L_{\overline i}
L_{\overline{i+1}} R_i)L_{\overline{i+1}}\dots L_{\overline k},
\label{pi}
\end{eqnarray}
\it
for all $k\geq 2$, $1\leq i\leq k-1$. Note, that all matrices entering these relations are of size $m^k\times m^k$.

Also, the degree $k\geq 2$ homogeneous component $\lr^{(k)}$ of the algebra $\lr$ is the quotient of $\LL^{\ot k}$ over the sum of the right hand side of (\ref{pi})
where $1\leq i \leq k-1$.
\end{proposition}

This construction is close to the transmutation procedure of \cite{M} and to that of \cite{IOP} used in the definition of quantum matrix algebras, associated with couples
of {\em compatible braidings}.

Since the maps $\pi_k$ are invertible, this proposition entails that any basis of the degree $k$ homogeneous component $\T^{(k)}$ of the algebra $\T$ is taken to a basis
of the component $\lr^{(k)}\subset \lr$.

\begin{remark} \label{re5}\rm
If the initial symmetry $R$ is involutive, it is so for the following braidings acting in the spaces $\TT^{\ot 2}$ and $\LL^{\ot 2}$:
\be
T_1\ot T_2\mapsto R^{-1}\, T_1\ot T_2\, R,\qquad L_{\overline 1}\ot L_{\overline 2}\mapsto R^{-1}\, L_{\overline{1}}\ot L_{\overline{2}}\,R.
\label{mapst}
\ee
By using this fact it is easy to construct
some symmetrizers  (projectors of symmetrisation) in the spaces $\TT^{(k)}$ and $\LL^{(k)}$ for any $k\ge 2$. Note that in the spaces $\LL^{(k)}$ the maps $\pi_k$ are involved in this construction.
This entails that the dimensions of the homogenous components $\T^{(k)}$ and $\lr^{(k)}$ for $k\ge 2$ are equal
to dimensions of the corresponding components ${\rm Sym}(gl(m))^{(k)}$ of the commutative algebra
${\rm Sym}(gl(m))$. Consequently, the algebras $\T$ and $\lr$ possess  the deformation property.

If an initial symmetry $R$ is Hecke, the  braidings (\ref{mapst})  acting in the spaces $\TT^{\ot 2}$ and $\LL^{\ot 2}$ have three eigenvalues, they are not symmetries at all. So,
the mentioned method of proving the deformation property should be modified. In \cite{GPS} we succeeded in constructing similar symmetrizers in the components $\TT^{(k)}$ and $\LL^{(k)}$
for $k=2,3$.
By applying these symmetrizers, it is possible to show that the dimensions of the homogeneous components $\T^{(k)}$ and $\lr^{(k)}$ for $k=2,3$
are equal to the dimensions of the corresponding components of the algebra
${\rm Sym}(gl(m))$ if $q-1$ is small enough. Though  similar symmetrizers in the higher components of these QM algebras are not still constructed,
according to \cite{D2} (see also \cite{PP}) it is possible to conclude  that for a generic $q$ the dimensions of the higher components $\T^{(k)}$, $\lr^{(k)}$, $k\ge 4$
are also equal to these in the algebra ${\rm Sym}(gl(m))$.
This entails the deformation property of the QM algebras, arising from Hecke symmetries.
\end{remark}

Now, describe the Poisson structures corresponding to the QM algebras. Below, we use the notation $\{L_1, L_2\}$ for the $m^2\times m^2$ matrix $L_1\, L_2$, where each entry
$l_i^j\ot l_k^l$ is replaced by $\{l_i^j, l_k^l\}$.

Let $R$ be a symmetry, which is a deformation of the flip $P$. Then the matrix $\RRR=PR$ is a deformation of the identity matrix:
\be \RRR=I-hr+O(h^2), \label{clas} \ee
where $h$ is a parameter of deformation and $r$ is the corresponding {\em classical $r$-matrix}, i.e. it obeys the relation
$$
[r_{12},r_{13}]+[r_{12},r_{23}]+[r_{13},r_{23}]=0.
$$
The negative sign in (\ref{clas}) is chosen for our convenience. If $R$ is a Hecke symmetry, we put $q=e^{-h}$.

If $R$ is an involutive (resp., Hecke) symmetry, then $r$ meets the following relation:
\be
r_{21}=-r_{12}\qquad ({\rm resp.,}\,\, r_{12}+r_{21}=2P).
\label{rrr}
\ee
In the latter case we present $r$ as a sum $r=r_-+r_+$, where
$$
r_- =\frac{r_{12}-r_{21}}{2},\qquad r_+=\frac{r_{12}+r_{21}}{2}
$$
are respectively the skew-symmetric and  symmetric components of $r$. Then the second relation in (\ref{rrr}) means that $r_+=P$. This relation is well known for the standard
Hecke symmetries, it is also valid for all other Hecke symmetries, for example, the Cremmer-Gervais ones.

The Poisson structure, corresponding to the RTT algebra, is
\be \{T_1, T_2\}=r T_1 T_2-T_2 T_1 r=r T_1 T_2-T_1 T_2 r. \label{Sk} \ee
Note that $T_1 T_2=T_2 T_1$ (and similarly for $L$) in virtue of the commutativity of the algebra ${\rm Sym}(gl(m))$ where the Poisson structure is defined.

If $r$ arises from an involutive symmetry, then it is clear that the right hand side of (\ref{Sk}) is skew-symmetric. If $r$ arises from a Hecke symmetry, then by using the relation $r_+=P$,
we arrive to the same conclusion. Moreover, $r$ in the bracket (\ref{Sk}) can be replaced by $r_-$.

Treating a classical $r$-matrix $r$ as an element of $gl(m)^{\ot 2}$, we can present the bracket (\ref{Sk}) as follows
\be
\{f,g\}=\cdot \big(\rho_l(r)^{\ot 2}(f\ot g) -\rho_r(r)^{\ot 2}(g\ot f)\big). \label{Skk}
\ee
 Here $\rho_l$ (resp., $\rho_r$) is the representation of $\gg$ by left (resp. right) vector fields acting onto the algebra ${\rm Sym}(gl(m))$
 and $\cdot$ stands for the commutative product in this algebra.

Now, consider the Poisson counterpart of the RE algebra $\lr$, associated with an involutive or Hecke symmetry $R$.
First, we rewrite the defining relations of the RE algebra $\lr$ in the following form
$$
\RRR_{12} L_1 \RRR_{21} L_2-L_2 \RRR_{12} L_1 \RRR_{21}=0.
$$

Applying the expansion (\ref{clas}) to the left hand side of this relation, we get the following Poisson bracket
\be
\{L_1, L_2\}=r_{12} L_1 L_2 +L_1 r_{21} L_2- L_2 r_{12} L_1 -L_2 L_1 r_{21}.
\label{STS}
\ee

If $r$ arises from an involutive symmetry, this bracket can be written as follows
\be
\{f,g\}=\cdot \rho_{ad}(r)^{\ot 2} (f\ot g),\quad f,g\in {\rm Sym}(gl(m)).
\label{STS1}
\ee
Here, $\rho_{ad}$ is the adjoint action of $gl(m)$ onto itself, extended upto the symmetric algebra ${\rm Sym}(gl(m))$ via the Leibniz rule. Thus, $\rho_{ad}^{\ot 2}(r)$ is
a skew-symmetric bi-vector field.

By contrast, if $r$ arises from a Hecke symmetry, a similar formula is as follows
\be
\{f,g\}=\cdot (\rho_{ad}^{\ot 2}(r_-)+(\rho_{r}\ot \rho_{l}-\rho_{l}\ot\rho_{r})(r_+)) (f\ot g).
\label{STS2}
\ee

It is clear that the bracket (\ref{STS2}) is skew-symmetric.
However, a straightforward verification of the Jacobi identity for it
is tedious. Nevertheless, the Jacobi identity is valid in virtue of the deformation property of the RE algebra.
Besides, the bracket (\ref{STS}) is compatible with the linear $gl(m)$ Poisson bracket, i.e. these two brackets form a Poisson pencil.
This claim can be easily verified by the lineatization procedure of the bracket (\ref{STS2}) (see Section 4).
Emphasize that the aforementioned algebra
$\lrh$ is the quantum counterpart of this pencil.

\begin{remark} \rm Let $\gg$ be a simple Lie algebra belonging to one of the classical series $B_m,\, C_m,$ $ D_m$, and $R$ be a braiding coming from the corresponding QG $U_q(\gg)$.
In this case the corresponding RTT and RE algebras are not deformations of the algebra ${\rm Sym}(\gg)$. The both brackets (\ref{Sk}) and (\ref{STS}) are
Poisson  on
 the corresponding Lie group $G$. Then they are respectively called the Sklyanin and the Semenov-Tian-Shansky brackets. Besides, they are also Poisson
on certain varieties (orbits) in $\gg^*$. We refer the reader to the paper \cite{DGS, Do}, where
these varieties are described.
\end{remark}

Now, we pass to the Yangians. Similarly to the construction above we introduce
the vector spaces
$$
\TT={\rm span}_{\K}(t_i^j[k]),\qquad \LL={\rm span}_{{\K}}(l_i^j[k])
$$
containing all finite linear combinations of the generators $t_i^j[k]$ of a Yangian of RTT type or, respectively, of the generators $l_i^j[k] $
of a braided Yangian, associated with a given
$R$-matrix (\ref{rrat}) or (\ref{rat}). Also, consider the spaces $\TT^{\ot k}$ and $\LL^{\ot k}$ containing all degree $k$ homogenous tensor polynomials
in the corresponding generators.

Define maps, similar to those in (\ref{pi}):
\be
\pi_k(T_1(u_1)\ot T_2(u_2)\ot ...\ot T_k(u_k))= L_{\overline 1}(u_1)\ot L_{\overline 2}(u_2)\ot ...\ot L_{\overline k}(u_k),
\label{pi1}
\ee
where as above we put
$$
L_{\overline 1}(u)=L_1(u),\qquad L_{\overline k}(u)=R_{k-1} L_{\overline {k-1}}(u) R^{-1}_{k-1}, \quad k\ge2.
$$
Here, we exhibited the maps  $\pi_k$ via the current matrices $T(u)$ and $L(u)$. Their presentation via the coefficients $T[s]$ and $L[s]$ has the same form.

Similarly to Proposition 4, the maps (\ref{pi1}) establish vector space isomorphisms between homogeneous components of the algebras $\YYR$ and $\YR$. As a consequence
we have the following.
\begin{proposition}\label{prop:7}
If a Yangian of RTT type has the deformation property, then the same is true for the corresponding braided Yangian.
\end{proposition}

In fact,  the method of the paper \cite{GPS} and mentioned in Remark 5 enables us to show the following claim.

\begin{proposition}
Any generalized Yangian is a deformation of the algebra  ${\rm Sym }(gl(m)[t^{-1}])$ provided that the initial symmetry
$R$ is a deformation of the flip.
\end{proposition}

A detailed proof of this proposition will be presented in our subsequent publication. Here we present some arguments in favor of this claim. In virtue of the Proposition \ref{prop:7}
we can restrict  ourselves to considering the Yangians of RTT type.

First,  observe that the map
\be
\rho: T_1(u)\ot T_2(v)\mapsto R(u,v)^{-1} T_1(v)\ot T_2(u)R(u,v) \label{mappp}
\ee
is involutive.  Indeed, we have $R(u,v)R(v,u)=\varphi(u,v) I$ where $\varphi(u,v)$ is a rational function depending on the initial symmetry $R$ (involutive or Hecke).
Thus, by applying the map $\rho$  twice we get the identity map.

This entails that the map $\rho$  has two eigenvalues $\pm 1$. The ideal coming in the definition
of the Yangian $\YYR$ is generated by the subspace $I_-=Im(I-R(u,v))$. Consider the complementary subspace $I_+=Im(I+R(u,v))$. Represent the defining relations of  this
subspace in the following form
$$R(u,v)T_1(u)\ot T_2(v)=-T_1(u)\ot T_2(v)R(u,v).$$

Expanding the generating matrix $T(u)$ in the usual way and using the relation
$$
\frac{1}{u-v}=\sum_{k\geq 0} v^k u^{-k-1},
$$
we get an explicit description of the space $I_+$ in terms of the entries the matrices
$T[k]$. Emphasize that this subspace is infinite dimensional and it is generated by quadratic polynomials in generators $t_i^j[r]$.

Consider the subspaces $I_{\pm}(h)$ in function of the deformation parameter $h$. As $h=0$ the subspaces $I_{\pm}(0)\subset\TT^{\ot 2}$ are respectively generated by
the elements
\be
t_j^c[r]  t_i^d[s] \pm  t_i^d[s] t_j^c[r].
\label{zero-h}
\ee
The generating series of these elements have the following matrix form
$$
T_1(u)\, T_2(v)\pm T_2(v)\, T_1(u).
$$

The subspace $I_+(h)$ is generated by the elements $E_{ij}^{cd}[k,l]$ which are coefficients at $u^{-k} v^{-l}$ of the folloving
generating series:
\be
E_{ij}^{cd}(u,v) = R_{ij}^{ab} (u,v)T_a^c(u) T_b^d(v)+T_i^a(v)T_j^b(u) R_{ab}^{cd}(u,v).
\label{lexi}
\ee
Note that at $h=0$ the elements $E_{ij}^{cd}[k,l]$ coincide with symmetric combinations of generators in (\ref{zero-h}). The subset of the elements $E_{ij}^{cd}[k,l]$ such that
the triple $(k, i, c)$ precedes that $(l,j, d)$ in the lexicographic order forms a basis of the space $I_+(0)$. Moreover, the space spanned by the elements $E_{ij}^{cd}[k,l]$
such that $k+l\leq p$ is finite dimensional.

Therefore, the elements $E_{ij}^{cd}[k,l]$ such that  $k+l\leq p$ are independent in the Yangian $\YYR$ if $h$  is small enough.  So, by sending the  elements $E_{ij}^{cd}[k,l]\in I_+(0)$ to their analogs in
$I_+(h)$ we construct the map $\alpha_h$ (see Introduction) at the quadratic level.
Reproducing this method we can construct the maps $\alpha_h$ for the higher homogenous components of the Yangian $\YYR$ and therefore show that this Yangian has the deformation property,
provided the symmetry $R$ is a deformation of the flip $P$.

Now compute the Poisson structures corresponding to the generalized Yangians.
We begin with the expansion $\RRR(u,v)=PR(u,v)=I-hr(u,v)+O(h^2)$ in a series in $h$, where $h=-a$ in formula (\ref{rrat}) and $h=-\log q$ in formula (\ref{rat}).
Then we get the following formulae for the classical current $r$-matrices
\be
r(u,v)=r-\frac{1}{u-v}P,\qquad r(u,v)=r-\frac{2u}{u-v}P, \label{ext}
\ee
where the constant matrix $r$ in the former formula corresponds to an involutive symmetry $R$ (and consequently, it is skew-symmetric), while $r$ in the latter formula
corresponds to a Hecke symmetry $R$ (and consequently, $r_{+}=P$).

Emphasize that the both current matrices $r(u,v)$ are skew-symmetric. Let us check this claim for the latter matrix:
$$
r_{12}(u,v)+r_{21}(v,u)=r_{12}+r_{21} - \frac{2u}{u-v} P -\frac{2v}{v-u}P=2r_+-2P=0.
$$
From now on, the notation $r_{21}(u,v)$ means that we interchange the elements from
 $gl(m)^{\ot 2}$ without interchanging $u$ and $v$.

The brackets, corresponding to the Yangians $\YR$ and $\YYR$ are  respectively similar to (\ref{Sk}) and (\ref{STS}):
\be \{T_1(u), T_2(v)\}=r(u,v)T_1(u) T_2(v)-T_2(v) T_1(u)r(u,v),   \label{SSk} \ee
\be \{L_1(u), L_2(v)\}=r_{12}(u,v)L_1(u) L_2(v)+L_1(u) r_{21} L_2(v)-L_2(v) r_{12} L_1(u)-L_2(v) L_1(u)r_{21}(u,v).  \label{SSTS} \ee

Observe that the middle $r$-matrices in (\ref{SSTS}) are constant, whereas the external factors are defined via one of the formulae  (\ref{ext}).

Since the Poisson structures under consideration are defined on the commutative algebra ${\rm Sym}(gl(m)[t^{-1}]$, we have
$$T_1(u)\, T_2(v)=T_2(v)\, T_1(u) ,\qquad L_1(u)\, L_2(v)=L_2(v)\, L_1(u).$$

\begin{remark} \rm
Emphasize that analogous brackets associated with Lie algebras of the series $B_m,$  $ C_m, $ $D_m$ are Poisson only on the corresponding groups.

Note that all brackets we are dealing with are  local, i.e. they  have no singularity
as $u-v \to 0$. More precisely, by putting $u-v=\nu$ and expanding $L(u)=L(v+\nu)$ in a series in $\nu$, we can present the above Poisson brackets via a series of brackets at
the same value $u=v$ but with the derivatives of the generating
matrices involved. At the quantum level this procedure is exhibited in \cite{GS}.
\end{remark}

\section{Linear Poisson brackets. Examples}

The linear Poisson brackets arising from the linearization procedure of the brackets above are also  of interest.
By applying this procedure we get the following linear counterpart of the bracket (\ref{SSk}) in the trigonometrical case
\be
\{T_1(u), T_2(v)\}=r(u,v)(T_1(u) +T_2(v))-(T_2(v) +T_1(u))r(u,v)=\Big[ r-\frac{2u}{u-v}P,\, T_1(u) +T_2(v)\Big ].
\label{ssk}
\ee
In the rational case we have
\be
\{T_1(u), T_2(v)\}=\Big[r-\frac{1}{u-v}\,P,\, T_1(u) +T_2(v)\Big].
\label{ssk1}
\ee

The linearization of the bracket (\ref{SSTS}) in the trigonometrical case leads to the bracket
\be
\{L_1(u), L_2(v)\}=-\Big[ 2P, \frac{u\,L_1(u) - v\,L_1(v)}{u-v}\Big].
\label{ssts}
\ee
In the rational case we have
\be
\{L_1(u), L_2(v)\}=-\Big[ P,\,\frac{L_1(u)-L_1(v)}{u-v} \Big].
\label{ssts1}
\ee

It is worth noticing that the bracket (\ref{ssts}) does not contain $r_-$ since the first term of the bracket (\ref{STS2}) vanishes under the linearization procedure. For the same
reason the constant term $r$ does not enter formula (\ref{ssts1}). Thus, these formulae do not depend of $r$, by contrast with the brackets (\ref{ssk}) and (\ref{ssk1}).

Note that the brackets (\ref{ssk1}) and (\ref{ssts1}) are similar. Namely, they coincide with each other if we put $r=0$ in (\ref{ssk1}) while the brackets (\ref{ssk}) and
 (\ref{ssts}) differ drastically.

The method of constructing the above linear brackets immediately entails that any of them is compatible with the corresponding quadratic Poisson bracket, i.e. these
two brackets (quadratic and linear) generate
a Poisson pencil. The quantum counterpart of such a Poisson pencil is the algebra which can be also obtained by the linearization of the corresponding Yangian.

Now, consider two examples.
First, we introduce the following symmetries:
\be
\hat R=\left(\begin{array}{cccc}
1&a&-a&ab\\
0&0&1&-b\\
0&1&0&b\\
0&0&0&1
\end{array}
\right), \qquad
R=\left(\begin{array}{cccc}
q&0&0&0\\
0&q-\qq&1&0\\
0&1&0&0\\
0&0&0&q
\end{array}
\right).
\label{Jord}
\ee
The former symmetry is involutive for any values of the parameters $a,b \in \K$. Note that it was introduced in \cite{G}.

The latter one is a Hecke symmetry. Note that it comes from the QG $U_q(sl(2))$.
The higher dimensional analogs of this Hecke symmetry come from the QG $U_q(sl(m))$, $m>2$. For all these Hecke symmetries the quantum determinants in
the corresponding RTT algebras are central.

By contrast, in the RTT algebra associated with the first symmetry in (\ref{Jord}) the determinant is central iff $a=b$. Consequently, the quantum determinant
in the corresponding Yangian of RTT type is also central iff $a=b$.

The current $R$-matrices corresponding to (\ref{Jord}) are of the form:
\be
\hat R(u,v)= \hat R-\frac{1}{u-v}I=
\left(\begin{array}{cccc}
1-\frac{1}{u-v}&a&-a&ab\\
0&-\frac{1}{u-v}&1&-b\\
0&1&-\frac{1}{u-v}&b\\
0&0&0&1-\frac{1}{u-v}
\end{array}\right).
\label{Rm}
\ee

\be
R(u,v)= R-\frac{(q-\qq)u}{u-v}I=
\left(\begin{array}{cccc}
\frac{-qv+\qq u}{u-v}&0&0&0\\
0&\frac{(-q+\qq)v}{u-v}&1&0\\
0&1&\frac{(-q+\qq)u}{u-v}&0\\
0&0&0&\frac{-qv+\qq u}{u-v}
\end{array}\right).
\label{Rm1}
\ee

The Yangian of RTT type, corresponding to the symmetry (\ref{Rm1}) is the lowest dimensional example of the aforementioned q-Yangians. A version of the
PBW bases in these Yangians can be extracted from the paper \cite{FJMR}. To this end it is necessary to consider a "half" of the quantum algebra  $U_q(\widehat{gl(m)})$ as defined in
\cite{FJMR}.

Now, consider the corresponding classical $r$-matrices. To this end, in the matrix $\hat R$ in (\ref{Jord}) we set $a=-h\alpha$, $b=-h\beta$, while for the matrix $R$ we set $q=e^{-h}$.
Expanding the $R$-matrices $\RRR=PR$ in $h$, we get the classical $r$-matrices:
\be
\begin{array}{l}
 \hat r=\al(e_{11}\otimes e_{12}-e_{12}\otimes e_{11})+\beta (e_{12}\otimes e_{22}-e_{22}\otimes e_{12})\\
\rule{0pt}{7mm}
r=e_{11}\ot e_{11}+e_{22}\ot e_{22} +2\,e_{21}\ot e_{12},
\end{array}
\label{exa}
\ee
where $e_{ij}$ are the elements of the usual basis of the Lie algebra $gl(2)$.

The $r$-matrix $\hat r$ in (\ref{exa}) is skew-symmetric. If $\al=\beta$ (that is $a=b$ in the symmetry $\hat R$ in (\ref{Jord})), it is an element of $sl(2)^{\wedge 2}$.
Unless, it is an element of $gl(2)^{\wedge 2}$. The second $r$-matrix is not skew-symmetric, its skew-symmetric component is $r_-=e_{21}\ot e_{12}-e_{12}\ot e_{21}$.

The classical $r$-matrices corresponding to the symmetries (\ref{Rm}) and (\ref{Rm1}) are respectively equal to
\be
\hat r(u,v)=\hat r -\frac{1}{u-v}P,\qquad r(u,v)=r-\frac{2u}{u-v}P.
\label{exa1} \ee

At last, we explicitly write down the linear Poisson brackets corresponding to the Yangian of RTT type associated with involutive symmetry $\hat R$ and to the braided Yangian
associated with the Hecke symmetry $R$. Denoting the entries of the matrix  $T(u)$ as follows
$$
T(u) = \left(\!\!
\begin{array}{cc}
a(u)&b(u)\\
\rule{0pt}{5mm}
c(u)&d(u)
\end{array}\!\!
\right),
$$
we get from (\ref{ssk1}):
$$
\begin{array}{lcl}
\{a(u),a(v)\} = \alpha(c(v)-c(u))&&
\{a(u),b(v)\} = \alpha(d(v)-a(v))+\frac{b(u)-b(v)}{u-v}\\
\rule{0pt}{6mm}
\{a(u),c(v)\} = \frac{c(v)-c(u)}{u-v}&&
\{a(u),d(v)\} =\beta c(u)-\alpha c(v)\\
 \rule{0pt}{7mm}
\{b(u),b(v)\} = (\alpha+\beta)(b(u)-b(v))&&
\{b(u),c(v)\} = (\alpha+\beta) c(v)+\frac{a(u)-a(v)+d(v)-d(u)}{u-v}\\
\rule{0pt}{6mm}
\{b(u),d(v)\} = \beta(d(u)-a(u))+\frac{b(u)-b(v)}{u-v}&&\\
\rule{0pt}{7mm}
\{c(u),c(v)\} = 0 && \{c(u),d(v)\} = \frac{c(v)-c(u)}{u-v}\\
\rule{0pt}{7mm}
\{d(u), d(v)\} = \beta (c(v)-c(u))
\end{array}
$$

Using the same notations for the entries of the generating matrix
$$
L(u) = \left(\!\!
\begin{array}{cc}
a(u)&b(u)\\
\rule{0pt}{5mm}
c(u)&d(u)
\end{array}\!\!
\right)
$$
of the braided Yangian $\YR$, we get from (\ref{ssts}):
$$
\begin{array}{lcl}
\{a(u),a(v)\} = 0&&
\{a(u),b(v)\} = \frac{2}{u-v}\,(ub(u)-vb(v))\\
\rule{0pt}{6mm}
\{a(u),c(v)\} = \frac{2}{u-v}\,(vc(v)-uc(u))&&
\{a(u),d(v)\} = 0\\
 \rule{0pt}{7mm}
\{b(u),b(v)\} = 0&&
\{b(u),c(v)\} = \frac{2}{u-v}\,(ua(u)-va(v)+vd(v)-ud(u))\\
\rule{0pt}{6mm}
\{b(u),d(v)\} = \frac{2}{u-v}\,(ub(u)-vb(v))&&\\
\rule{0pt}{7mm}
\{c(u),c(v)\} = 0 && \{c(u),d(v)\} = \frac{2}{u-v}\,(vc(v) - uc(u))\\
\rule{0pt}{7mm}
\{d(u), d(v)\} = 0
\end{array}
$$
Note that the Poisson structure (\ref{ssts1}) of the braided Yangian for the rational current $R$-matrix $\hat R(u,v)$ can be easily obtained from the above formulae for $T(u)$
just by setting $\alpha =0$ and $\beta = 0$.
\begin{remark}\rm
 Note that the construction of the affine quantum algebra $U_q(\widehat{gl(m)})$ from \cite{FJMR} can be directly generalized to other current $R$-matrices.
 Let us precise that this construction is based on matching two Yangian-type algebras via
 some {\em permutation relations}. A similar method can be applied in order to construct a "braided version" of this algebra from two braided Yangians.
 However, it is not clear what are centers of these algebras.

 \end{remark}

\end{document}